\documentclass{svjour3}                     
\smartqed  
\usepackage{graphicx}

\begin{document}

\title{A harmonic Lanczos bidiagonalization method for computing interior singular triplets of large
matrices}

\titlerunning{harmonic method for interior SVD problems}        

\author{Datian Niu \and
        Xuegang Yuan}


\institute{D. Niu \at
              School of Science, Dalian Nationalities University, \\
              Dalian, 116600, China\\
              \email{niudt@dlnu.edu.cn}           
           \and
           X. Yuan \at
              School of Science, Dalian Nationalities University, \\
              Dalian, 116600, China\\
              \email{yuanxg@dlnu.edu.cn}           
           \and
              This research was Supported by the NSFC Grants
              10872045 and by Program for New Century Excellent Talents in University.
           }

\date{Received: date / Accepted: date}

\maketitle

\begin{abstract}
This paper proposes a harmonic Lanczos bidiagonalization method for
computing some interior singular triplets of large matrices. It is shown
that the approximate singular triplets are convergent if a certain Rayleigh quotient matrix
is uniformly bounded and the approximate
singular values are well separated. Combining with the implicit
restarting technique, we develop an implicitly restarted harmonic
Lanczos bidiagonalization algorithm and suggest a selection strategy of shifts.
Numerical experiments show that one can use this algorithm to compute interior
singular triplets efficiently.

\end{abstract}

\section{Introduction}
\label{intro}

The singular value decomposition (SVD) of a matrix $A\in R^{M\times
N}, M\geq N$ is given by
\begin{equation}
A=U\Sigma V^{\rm T},
\end{equation}
where $\Sigma=diag(\sigma_1, \sigma_2, \cdots, \sigma_N)$,
$U=(u_1,u_2,\cdots,u_M)$ and $V=(v_1,v_2,\cdots,v_N)$ are orthogonal
matrices of order $M$ and $N$ respectively. $(\sigma_i,u_i,v_i),i=1,2,\cdots N,$ are
called the singular triplets of $A$.

Consider the $(M+N)\times (M+N)$ augmented matrix
\begin{equation}
\tilde A = \left ( \begin{array}{cc}
                     0 & A \\
                     A^{\rm T} & 0
                     \end{array} \right ).\label{aug}
\end{equation}
Then, the eigenvalues of $\tilde A$ are $\pm \sigma_1, \pm \sigma_2,
\cdots, \pm \sigma_N$ and $M-N$ zeros. The eigenvectors associated
with $\sigma_i$ and $-\sigma_i$ are $\frac{1}{\sqrt{2}}\left
(u_i^{\rm T},v_i^{\rm T}\right )^{\rm T}$ and
$\frac{1}{\sqrt{2}}\left (u_i^{\rm T},-v_i^{\rm T}\right )^{\rm T}$
respectively. Therefore, the SVD problems are equivalent to the
eigenproblems of augmented matrices.

The SVD methods are widely used in determination of numerical rank,
determination of spectral condition number, least square problems,
regression analysis, image processing, signal processing, pattern
recognition, information retrieval, and so on.

At present, computation of largest or smallest singular triplets of
large matrices has been well studied, Lanczos bidiagonalization
method and its variants are the most popular methods.
In 1981, Golub et al. \cite{golub1981} firstly designed a block
Lanczos bidiagonalization method to compute some largest singular triplets.
Larsen \cite{larsen1998} discussed the reorthogonalization of the Lanczos
bidiagonalization process. Jia and Niu \cite{jianiu2003} proposed a
refined Lanczos bidiagonalization method to compute some largest and
smallest singular triplets. Kokiopoulou et al. \cite{kokio2004} used the harmonic
projection technique to compute the smallest singular values. Baglama and Reichel
\cite{baglama2005,baglama2006} used Ritz values and harmonic Ritz values to
approximate the largest and smallest singular values respectively.
Hernandez et al. \cite{hernandez2008} provided a parallel implementation of the Lanczos bidiagonalization method.
Stoll \cite{stoll2008} developed a Krylov-schur approch to partial SVD.
Recently, Jia and Niu \cite{jianiu2009} proposed a refined harmonic Lanczos bidiagonalization method
to compute some smallest singular triplets.
All of above methods compute the Lanczos bidiagonalization process,
build two $m-$dimensional Krylov subspaces, then extract approximate
singular triplets from these two subspace by different ways.
Hochstenbach \cite{hochstenbach2001,hochstenbach2004} also give
the Jacobi-Davidson type algorithms for SVD problems.

Due to the storage requirement and the computational cost, all the
projection methods must be restarted. The implicit restarting
technique \cite{sorensen1992} proposed by Sorensen is the most
powerful tool and is widely used in many projection methods. The
success of this technique heavily depends on the selection of the shifts, see \cite{jia1999,sorensen1992}.
For eigenvalue problems, Sorensen \cite{sorensen1992} used the unwanted
Ritz values as the shifts to restart Arnoldi method, and Morgan
\cite{morgan2006} used the unwanted harmonic Ritz values as the shifts to
restart harmonic Arnoldi method. Jia \cite{jia1999,jia2002} used
the refined shifts and refined harmonic shifts obtained by the
information of the refined Ritz vectors and refined harmonic vectors to
restart refined Arnoldi method and refined harmonic Arnoldi method,
respectively. For SVD problems, Kokiopoulou et al. \cite{kokio2004}
used the unwanted harmonic Ritz values as the shifts. Baglama and Reichel
\cite{baglama2005,baglama2006} explicitly augmented the Lanczos bidiagonalization method
with certain Ritz vectors or harmonic Ritz vectors. Jia and Niu \cite{jianiu2003,jianiu2009}
gave an refined (harmonic) shift strategy within the implicitly restarted refined (harmonic)
Lanczos bidiagonalization method.

In this paper, we are concerned with the computation of interior
singular triplets. For a given target $\tau$, we want to compute
some singular triplets nearest $\tau$. So, we sort the singular
triplets by
\begin{equation}\label{sort}
|\sigma_1-\tau|\le |\sigma_2-\tau| \le \cdots \le |\sigma_N-\tau|.
\end{equation}

We must emphasize that, in this paper, $\sigma_1$ is the singular value nearest $\tau$
rather than the smallest singular value, meanwhile, $\sigma_N$ is the singular value farthest
from $\tau$ rather than the largest singular value.

Since the largest eigenvalues of $(\tilde A-\tau I)^{-1}$ are the
eigenvalues of $\tilde A$ closest to $\tau$, and the SVD problem of
$A$ is equivalent to the eigenproblem of $\tilde A$, we can use
shift-invert technique on $\tilde A-\tau I$ to compute the interior
singular triplets, such as shift-and-invert Arnoldi method ({\sf svds}).
In this paper, we assume that $M$ and $N$ are large and that $A$ can not be factorized.
The shift-and-invert technique need the factorization of $\tilde A-\tau I$.
Since $M$ and $N$ are large, $M+N$, the dimension of $\tilde A-\tau I$, is larger.
We can not do any factorizations on $\tilde A-\tau I$. Therefore, the shift-and-invert technique is not
suitable for interior SVD problems.

Another approach for computing interior singular triplets is the harmonic
projection method. The harmonic projection method has been widely used
to compute interior eigenpairs, see \cite{morgan2000,morgan2006},
and has been adopted to combine with Lanczos bidiagonalization methods to compute smallest singular triplets \cite{baglama2005,baglama2006,kokio2004,jianiu2003}.
However, if we use the harmonic projection method explicitly on $\tilde A-\tau I$,
the scale of the problem is increased and this leads to the increasing computational cost.
Further, we ignore the special structure of $\tilde A$ or $\tilde A-\tau I$, and the projected matrix
and the updated process of implicit restarting may lose this
structure. Therefore, we must use the harmonic projection method implicitly.
Until now, no literature has been appeared to compute interior singular
triplets by the harmonic projection method implicitly.

In this paper, we propose a harmonic Lanczos bidiagonalization method for computing interior singular triplets 
by combining the harmonic projection technique with the Lanczos bidiagonalization process. 
We analyze the convergence behavior,
show that the harmonic Ritz approximations converge to the desired interior singular triplets if
some Rayleigh quotient matrix is uniformly bounded and the harmonic Ritz values are
well separated. Then, based on Morgan's harmonic shift strategy
\cite{morgan2006} for computing interior eigenvalues, we give a
selection of the shifts within the framework of the implicitly restarted
harmonic Lanczos bidiagonalization methods. Further, we report some
numerical experiments of computation of interior singular triplets.
It appears that the algorithm we proposed is suitable for computing
the interior singular triplets of large matrices.

Throughout this paper, denote by $||\cdot||$ the spectral norm of a
matrix and the vector 2-norm, by ${\cal
K}_m(C,v_1)=span\{v_1,Cv_1,\cdots,C^{m-1}v_1\}$ the $m-$ dimensional
Krylov subspace generated by the matrix $C$ and the starting vector
$v_1$, by superscript '$^{\rm T}$' the transpose of matrix or
vector, by $e_m$ the $m-$th coordinate vector of dimension $m$.

\section{Harmonic Lanczos bidiagonalization method}

\subsection{Lanczos bidiagonalization process}

Golub et al. \cite{golub1981} proposed a Lanczos bidiagonalization
method to compute the largest singular triplets of $A$. This method
is equivalent to the symmetric Lanczos method on $\tilde A$ with a
special initial vector. It is based on the Lanczos bidiagonalization
process, which is shown in matrix form as follows:
\begin{eqnarray}
AQ_m=P_mB_m, \label{bid1}\\
A^{\rm T}P_m=Q_mB_m^{\rm T}+\beta_m q_{m+1}e_m^{\rm T}, \label{bid2}
\end{eqnarray}
where
\begin{equation}
    B_m=\left ( \begin{array} {c c c c}
                       \alpha_1 & \beta_1  &        &         \\
                                & \alpha_2 & \ddots &         \\
                                &          & \ddots & \beta_{m-1} \\
                                &          &        & \alpha_m
                     \end{array} \right)
\end{equation}
is an upper bidiagonal matrix, $Q_m=(q_1,q_2,\ldots,q_m)$ and
$P_m=(p_1,p_2,\ldots,p_m)$ span the Krylov subspaces ${\cal
K}_m(A^{\rm T}A,q_1)$ and ${\cal K}_m(AA^{\rm T},p_1)$,
respectively.

In finite precision arithmetic, the columns of $P_m$ and $Q_m$ may
lose the orthogonality rapidly and must be reorthognalized. From
the analysis of Simon and Zha \cite{simonzha2000}, 
we know that only the columns of one of the matrices $P_m$ and $Q_m$ need to be reorthogonalized.
When $M\gg N$, Reorthogonalization on $Q_m$ only can reduce the
computational cost considerably. So we only perform
reorthogonalization on $Q_m$.

\subsection{Harmonic Lanczos bidiagonalization method}
Given the subspace
\begin{equation}
{\cal E}=span\left \{ \left ( \begin{array}{cc}
                         P_m & 0\\
                         0 & Q_m
                       \end{array} \right ) \right \}.
\end{equation}
Making use of the harmonic projection principle, we compute some
approximate eigenpairs $(\theta_i,\tilde \varphi_i)$ of $\tilde A$
nearest $\tau$ by requiring
\begin{equation}
\left \{ \begin{array}{c}
         \tilde \varphi_i\in {\cal E},\\
         (\tilde A-\theta_i I)\tilde \varphi_i \bot (\tilde A - \tau
         I){\cal E}.
         \end{array} \right . \label{harmproj}
\end{equation}

From (\ref{bid1}) and (\ref{bid2}), (\ref{harmproj}) can be
rewritten as the following generalized eigenproblem:
\begin{equation}
\left ( \begin{array}{cc}
           -\tau I & B_m \\
           B_m^{\rm T} & -\tau I
        \end{array} \right )
\left ( x_i \atop y_i \right )=\frac{1}{\theta_i-\tau} \left (
\begin{array}{cc}
           \tau^2 I+B_mB_m^{\rm T}+\beta_me_me_m^{\rm T} & -2\tau B_m\\
           -2\tau B_m^{\rm T} & \tau^2 I+B_m^{\rm T}B_m
        \end{array} \right )
\left ( x_i \atop y_i \right ).
\end{equation}

Assume that $\theta_i>0, i=1,2,\cdots,k+l$, which are sorted by
$$|\theta_1-\tau | \le
|\theta_2-\tau | \le \cdots \le |\theta_{k+l}-\tau |$$ and
$\theta_i<0,i=k+l+1,k+l+2,\cdots,2m$. We can use
$\theta_i,i=1,2,\cdots,k$ and $\tilde \varphi_i=\left (P_mx_i \atop
Q_my_i \right )$ as the approximation of the desired eigenpair of
$\tilde A$. Because of the relation between the singular triplets of $A$
and the eigenpairs of $\tilde A$, we use $\theta_i, \tilde
u_i=P_mx_i/||x_i||=P_m\tilde x_i, \tilde
v_i=Q_my_i/||y_i||=Q_m\tilde y_i,i=1,2,\cdots,k$ as the approximate
singular triplets of $A$ nearest $\tau$. Here we call $\theta_i, \tilde
u_i, \tilde v_i$ the harmonic Ritz value, the left and right harmonic Ritz
vector, respectively.

From (\ref{bid1}) and (\ref{bid2}), we have
$$||A\tilde v_i-\theta_i\tilde u_i||=||B_m\tilde y_i-\theta_i\tilde
x_i||,$$
$$ ||A\tilde v_i-\theta_i\tilde u_i||=\sqrt{||B_m^{\rm T}\tilde
x_i-\theta_i\tilde y_i||^2+\beta_m^2|e_m^{\rm T}\tilde x_i|^2}.$$

Therefore, if
\begin{equation}\label{stop}
\sqrt{||B_m\tilde y_i-\theta_i\tilde x_i||^2+||B_m^{\rm T}\tilde
x_i-\theta_i\tilde y_i||^2+\beta_m^2|e_m^{\rm T}\tilde x_i|^2}<tol,
\end{equation}
where $tol$ is a prescribed tolerance, then the method is known as
convergent. So we need not form $\tilde u_i$ and $\tilde v_i$
explicitly before convergence.

\subsection{Convergence analysis}

Set
$$\tilde B=\left ( \begin{array}{cc}
                        -\tau I & B_m \\
                        B_m^{\rm T} & -\tau I
                     \end{array} \right )$$
and
$$\tilde C=\left ( \begin{array}{cc}
           \tau^2 I+B_mB_m^{\rm T}+\beta_me_me_m^{\rm T} & -2\tau B_m\\
           -2\tau B_m^{\rm T} & \tau^2 I+B_m^{\rm T}B_m
        \end{array} \right ),$$
then $\theta_i, i=1, 2, \cdots, 2m$ are the eigenvalues of $\tilde
B^{-1}\tilde C$. The matrix $\tilde B$ is called the Rayleigh quotient matrix of $\tilde A$
with respect to the subspace ${\cal E}$ and the target $\tau$.

The following results are direct from Theorem 2.1, Corollary 2.2 and
Theorem 3.2 of \cite{jia2005}.

\begin{theorem}\label{Thvalue}
Assume that $(\sigma,u,v)$ is a singular triplet of $A$, define
that $\epsilon=\sin\angle\left (\left (u\atop v \right ),{\cal E}\right )$
is the distance between the vector $\left (u\atop v \right )$ and
the subspace ${\cal E}$. Then there exists a perturbation matrix F such
that $\sigma$ is an exact eigenvalue of $\tilde B^{-1}\tilde C+F$,
where
\begin{equation}
||F||\le\frac{\epsilon}{\sqrt{1-\epsilon^2}}||\tilde
B^{-1}||(\sigma||A||+||A||^2).
\end{equation}
Furthermore, there exists an eigenvalue of $\tilde B^{-1}\tilde C$
satisfying
\begin{equation}
|\theta-\sigma|\leq(2||A||+||F||)||F||.
\end{equation}
\end{theorem}

Theorem \ref{Thvalue} shows that if $\epsilon$ tends to zero and if
$||\tilde B^{-1}||$ is uniformly bounded, then there exists one
harmonic Ritz value $\theta$ converging to the desired singular
value $\sigma$.

However, from the interlacing theorem of eigenvalues
\cite{golub1996}, since
$$\tilde B=\left ( \begin{array}{cc}
                         P_m & 0\\
                         0 & Q_m
                       \end{array} \right )^{\rm T}
           (\tilde A-\tau I)
           \left ( \begin{array}{cc}
                         P_m & 0\\
                         0 & Q_m
                       \end{array} \right ),$$
we have that the eigenvalues of $\tilde B$ are between the largest and smallest eigenvalue of
$\tilde A -\tau I$. Therefore, $\tilde B$ may be singular, which leads to arbitrarily large
$\|\tilde B^{-1}\|$. Hence, we must assume $\|\tilde B^{-1}\|$ is uniformly bounded.
In fact, this is the inherent defect of the harmonic projection methods, which can be easily obtained
from Jia's analysis \cite{jia2005}.

Similarly to the analysis in \cite{jia2005}, if $\tau$ is very close
to a desired singular value $\sigma$ of $A$, then the method may
miss it. We replace $\theta_i$ by the Rayleigh-quotient
$\rho_i=\tilde u_i^TA\tilde v_i=\tilde x_i^TB_m\tilde y_i$ as the
approximate singular value, as was done in
\cite{hochstenbach2004,jianiu2009}. In general, $\rho_i$ is more
accurate than $\theta_i$.

\begin{theorem}\label{Thvector}
Let $(\theta,z)$ be an eigenpair of $\tilde B^{-1}\tilde C$, where
$z=\left (x \atop y\right )$, and assume $(z,Z_\bot)$ to be
orthogonal such that
\begin{equation}
\left (z^{\rm T} \atop Z_\bot^{\rm T} \right )\tilde B^{-1}\tilde
C(z,Z_\bot)= \left (\begin{array}{cc}
          \theta & g^{\rm T} \\
          0 & G
       \end{array} \right ).
\end{equation}
If
\begin{equation}
sep(\theta,G)=||(G-\theta I)^{-1}||^{-1}>0,
\end{equation}
then
\begin{eqnarray}
\sin \angle \left (\left (u \atop v \right ),\left ( \tilde u \atop
\tilde v \right )\right ) & \le &\left (1+\frac{2||\tilde
B^{-1}||||A||}{\sqrt{1-\epsilon^2}sep(\sigma,G)}
\right )\varepsilon \nonumber \\ 
&\le & \left (1+\frac{2||\tilde
B^{-1}||||A||}{\sqrt{1-\epsilon^2}(sep(\theta,G)-|\sigma-\theta|)}
\right )\varepsilon.
\end{eqnarray}
\end{theorem}

Theorem \ref{Thvector} shows that if $\|\tilde B^{-1}\|$ is uniformly bounded and $sep(\theta,G)$ is bounded
below by a positive constant, that is, all harmonic Ritz values are
well separated, then the harmonic Ritz vectors $\tilde
u,\tilde v$ converge to the desired left and right singular vector.

\section{Implicit restarting technique, shifts selection and an adaptive shifting strategy}

\subsection{Implicit restarting technique}

Due to the storage requirement and the computational cost, the number of
Lanczos bidiagonalization steps $m$ can not be large. However, for a
relatively small $m$, the approximate singular triplets do not converge.
Therefore, the method must be restarted generally.

The implicit restarting technique proposed by Sorensen
\cite{sorensen1992} is a powerful restarting tool for the Lanczos and
Arnoldi process, and has been adopted to the Lanczos bidiagonalization
process \cite{bjorck1994,jianiu2003,jianiu2009,kokio2004,larsen}.
After running the implicit QR iteration $p$ steps on $B_m$ and
using the shifts $\mu_j,j=1,2,\cdots,p$, we have
\begin{equation}
\left \{ \begin{array}{l}
           (B_m^{\rm T}B_m-\mu_1^2I)\cdots(B_m^{\rm T}B_m-\mu_p^2I)=\tilde P
           R,\\
           \tilde P^{\rm T} B_m \tilde Q  \mbox{ upper bidiagonal},
         \end{array} \right .
\end{equation}
where $\tilde P, \tilde Q$ are the products of the left and right
Givens rotation matrices applied to $B_m$.

Performing the above process gives the following relation:
\begin{eqnarray}
AQ_{m-p}^+&=&P_{m-p}^+B_{m-p}^+,\\
A^{\rm T}P_{m-p}^+&=&Q_{m-p}^+{B_{m-p}^+}^{\rm  T}+(\beta_{m-p}\tilde p_{m,m-p}
q_{m+1}+\beta_{m-p}^+ q_{m-p+1}^+)e_{m-p}^{\rm  T},
\end{eqnarray}
where $Q_{m-p}^+$ and $q_{m-p+1}^+$ are the first $m-p$ columns and the $(m-p+1)$-th
column of $Q_m\tilde Q$, $P_{m-p}^+$ is the first $m-p$ columns of
$P_m\tilde P$, $B_{m-p}^+$ is the leading $(m-p)\times (m-p)$ block of $\tilde P
B_m \tilde Q$, $\tilde p_{m,m-p}$ is the $(m,m-p)$ element of $\tilde
P$. Since $\beta_{m-p}\tilde p_{m,m-p}q_{m+1}+\beta_{m-p}^+ q_{m-p+1}^+$ is
orthogonal to $Q_{m-p}^+$, we obtain a $(m-p)$-step Lanczos
bidiagonalization process starting with $q_1^+$, where
\begin{equation}\label{initq}
\gamma q_1^+=\prod_{j=1}^p(A^{\rm T}A-\mu_j^2 I)q_1 \label{update}
\end{equation}
with $\gamma$ a factor making $\|q_1^+\|=1$. It is then extended to
the $m$-step Lanczos bidiagonalization process in a standard way.

\subsection{shifts selection and adaptive shifting strategy}

Once the shifts $\mu_1,\mu_2,\ldots,\mu_p$ are given, we can run the
implicitly restarted algorithm described above iteratively. The success of
the implicit restarting technique heavily depends on the selection
of the shifts. As is shown in \cite{jianiu2003}, from (\ref{initq}), it
can be easily seen that the more accurate the shifts approximate to
some unwanted singular values, the more information on the unwanted
singular vectors are dampened out after restarting. Therefore, the
resulting subspace contains more information on the desired singular
vectors, and the algorithms may converge faster. For eigenproblems
and SVD problems, Morgan \cite{morgan2006} and Kokiopoulou et al.
\cite{kokio2004} suggested using the unwanted harmonic Ritz values
as shifts. A natural choice of the shifts within our algorithm is
the unwanted approximate singular values
$\theta_{k+j},j=1,2,\cdots,l$, since they are the best approximations
available to some of the unwanted singular values within our
framework.

From (\ref{initq}), we see the component along the desired $k$-th
singular vector $u_k$ is greatly damped if a shift $\mu_i$ is very
close to $\sigma_k$, so $\mu_i$ is a bad shift and $\rho_k$
may converge to $\sigma_k$ very slowly or not at all. To correct this
problem, Larsen \cite{larsen} proposed an adaptive strategy to
compute largest singular triplets. He replaces a bad shift by
zero shift. Jia and Niu \cite{jianiu2003,jianiu2009} gave a modified form
for computing smallest singular triplets. Define the relative gaps
of $\rho_k$ and all the shifts $\mu_i,i=1,2,\cdots,l$ by
  \begin{equation}
    {\rm  relgap}_{ki}=\left |\frac{(\rho_k-\varepsilon_k)-\mu_i}{\rho_k}\right |,
    \label{harmrelgap}
    \end{equation}
where $\varepsilon_k$ is the residual norm (\ref{stop}). If ${\rm
relgap}_{ki}\leq 10^{-3}$, $\mu_i$ is a bad shift and should be
replaced by a suitable quantity. They replace the bad shifts by the
largest or the smallest approximate singular value for computing the
smallest or the largest singular triplets. In this paper, a good
strategy is replacing the bad shifts by the approximate
singular value farthest from $\tau$, as this strategy amplifies
the components of $q_1^+$ in $v_i,i=1,2,\cdots,k$ and damps those in
$v_i,i=k+1,k+2,\cdots,N$.

\section{Numerical Experiments}

Numerical experiments are carried out using Matlab 7.1 R14 on an
Intel Core 2 E6320 with CPU 1.86GHZ and 2GB of memory under the
Window XP operating system. Machine epsilon is $\epsilon_{\rm
mach}\approx 2.22\times 10^{-16}$. The stopping criteria is
\begin{equation}
   stopcrit=\max_{1 \leq i \leq k}\sqrt{\|A\tilde v_i-\rho_i
    \tilde u_i\|^2+\|A^{\rm T}\tilde u_i-\rho_i\tilde v_i\|^2}.
\end{equation}
If
\begin{equation}
\frac{stopcrit} {\|A\|_1}<tol, \label{stopcrit}
\end{equation}
then stop. From (\ref{stop}), we need not form $\tilde u_i, \tilde
v_i$ explicitly before convergence.

For large eigenproblems, in order to speed up convergence, most of the
implicitly restarted Krylov type subspace algorithms, such as
ARPACK({\sf eigs}), compute $k+3$ approximate eigenpairs when $k$
eigenpairs are desired. This strategy has been adopted to SVD problems,
see \cite{baglama2005,jianiu2009}. In this paper, we also compute $k+3$
approximate singular triplets and use $l-3$ shifts in implicit
restarting process.

All test matrices are from \cite{bai1997}. We take $tol=10^{-6}$. In
all the tables, '$iter$' denotes the number of restart, '$time$'
denotes the CPU timings in second, '$mv$' denotes the number of
matrix-vector products. Since the matrix-vector products performed on
$A$ are equal to those on $A^{\rm T}$, we only count the matrix-vector
products on $A$.

\subsection{Computation of smallest singular triplets}

Obviously, we can compute some smallest singular triplets by taking $\tau=0$.
We compute three singular triplets nearest $\tau=0, 0.01, 0.005, 0.001$
of WELL1850, respectively. These three singular values are all the
three smallest singular values. The computed three singular values
are
$$\sigma_1 \approx 1.611969e-002, \sigma_2 \approx 1.911309e-002,
\sigma_3 \approx 2.315889e-002.$$

\begin{figure*}
  \includegraphics[width=0.5\textwidth]{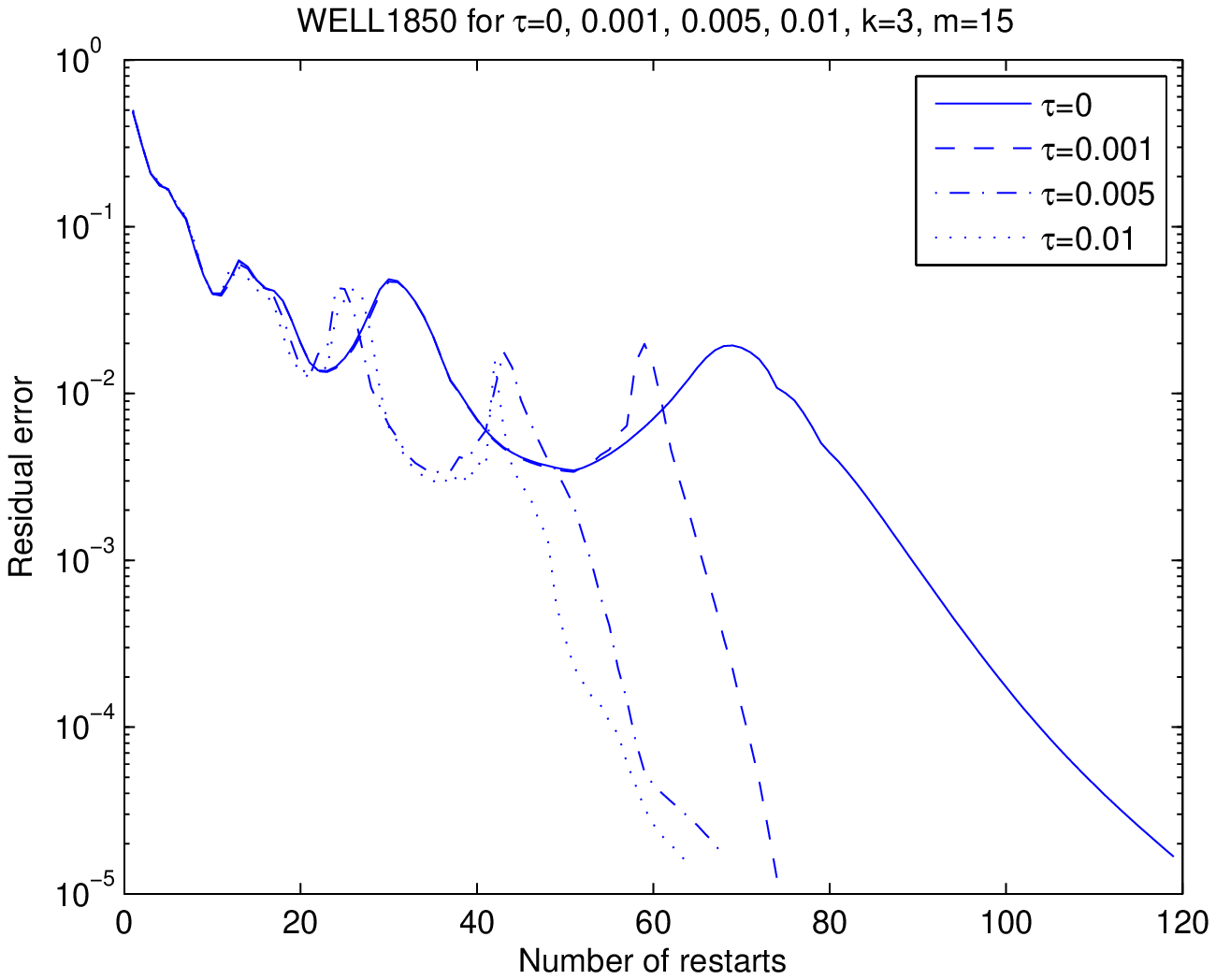}
  \includegraphics[width=0.5\textwidth]{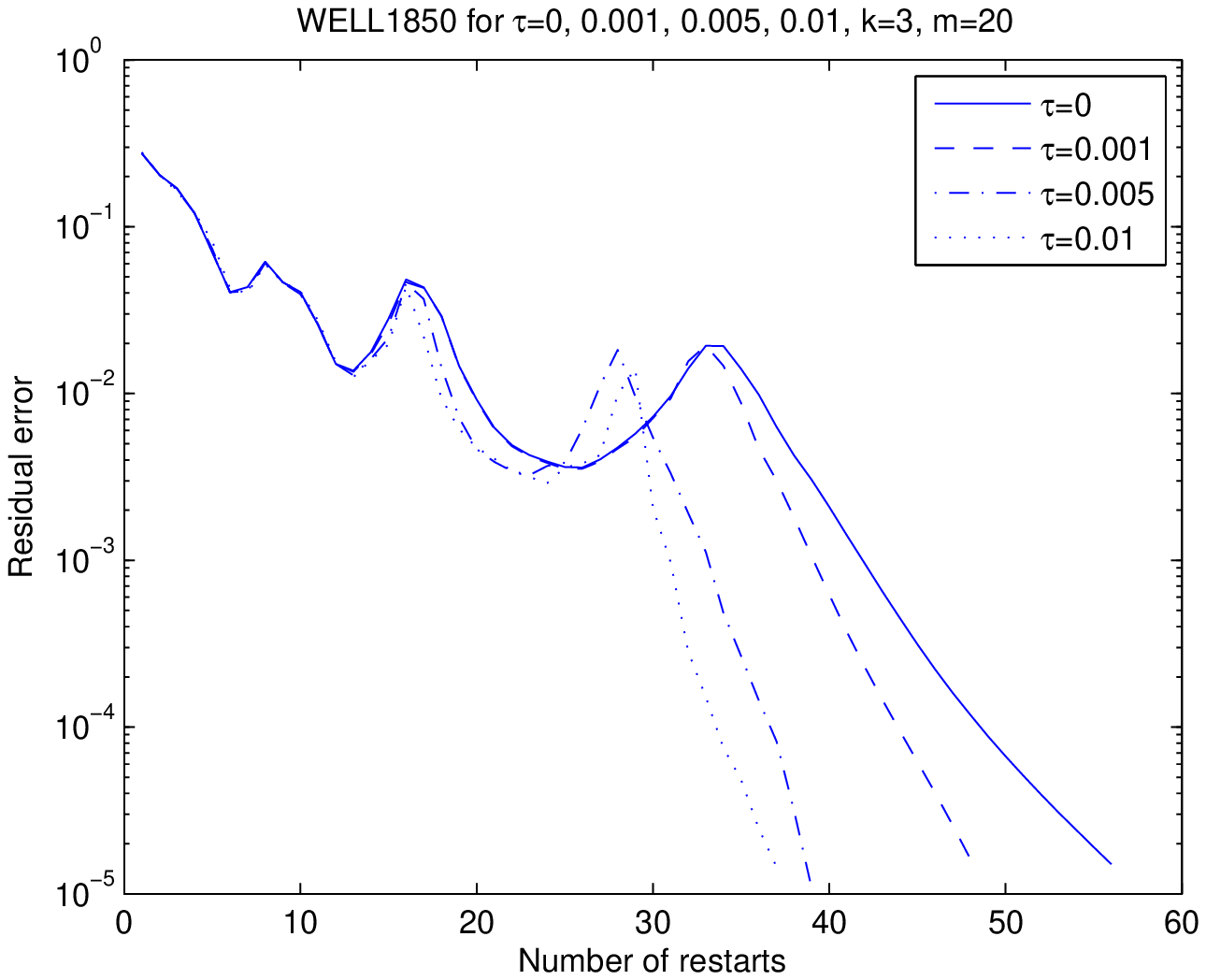}
  \caption{Absolute residual norms for WELL1850 for $k=3$, $m=20$, $\tau=0,0.01,0.005,0.001$}
  \label{fig1}
\end{figure*}

Table \ref{tab1} reports the computational results. Fig. \ref{fig1}
plots the absolute residual norms of the computed singular triplets for
$m=15$ and $m=20$, respectively. From Table \ref{tab1} and Fig.
\ref{fig1}, we see that for all $\tau$, our algorithm can compute
three singular triplets accurately. However, for different $\tau$,
the algorithm has a great difference on restart numbers,
matrix-vector products and CPU times. This phenomenon shows a good
choice of target point $\tau$ can speed up the convergence
considerably.

\begin{table}[htp]
\caption{WELL1850 for $k=3$, $m=10,15,20,25$, $\tau=0,0.001,0.005,0.01$}
\label{tab1}
\begin{tabular}{|c|c|c|c|c|c|c|c|c|} \hline
$m$ & $iter$ & $time$ & $mv$ & $stopcrit$ & $iter$ & $time$ & $mv$ & $stopcrit$ \\ \hline
&\multicolumn{4} {|c|} {$\tau=0$}&\multicolumn{4} {|c|} {$\tau=0.001$} \\ \hline
 10 & 543 & 5.15 & 2178 & 1.67e-005 & 178 & 1.89 & 718 & 1.66e-005 \\ \hline
 15 & 119 & 3.20 & 1077 & 1.67e-005 & 74 & 1.99 & 672 & 1.25e-005  \\ \hline
 20 & 56 & 3.48 & 790 & 1.50e-005 & 48 & 2.81 & 678 & 1.62e-005  \\ \hline
 25 & 35 & 3.20 & 671 & 1.35e-005 & 35 & 3.64 & 671 & 1.14e-005  \\ \hline
&\multicolumn{4} {|c|} {$\tau=0.005$}&\multicolumn{4} {|c|} {$\tau=0.01$} \\ \hline
 10 & 160 & 1.66 & 646 & 1.66e-005 & 179 & 1.76 & 722 & 1.60e-005 \\ \hline
 15 & 68 & 1.86 & 618 & 1.62e-005 & 64 & 1.73 & 582 & 1.52e-005  \\ \hline
 20 & 39 & 2.31 & 552 & 1.08e-005 & 37 & 2.23 & 524 & 1.45e-005  \\ \hline
 25 & 31 & 3.10 & 595 & 1.27e-005 & 29 & 2.89 & 557 & 1.06e-005  \\ \hline
\end{tabular}
\end{table}

\subsection{Computation of three interior singular triplets nearest different $\tau$}

The test matrix is DW2048, a $2048 \times 2048$ matrix. We compute
three singular triplets nearest different $\tau$. The
computational results are shown in Tables \ref{tab2}-\ref{tab3}.
From Table \ref{tab2}, we see that the relative errors of the computed
singular values are no more than $O(10^{-9})$. The Tables
demonstrate that our algorithm can compute the desired singular
triplets accurately.

\begin{table}[htp]
\caption{Three computed singular values of DW2048 nearest $\tau=0.2,
0.5, 0.6, 0.8$ for $m=50$}\label{tab2}
\begin{tabular}{|c|c|c|c|} \hline
 \multicolumn{2} {|c|} {$\tau=0.2$} & \multicolumn{2} {|c|} {$\tau=0.5$}
 \\ \hline
 $\rho_j$ & $|\rho_j-\sigma_j|/\sigma_j$ & $\rho_j$ & $|\rho_j-\sigma_j|/\sigma_j$ \\ \hline
 2.0031301e-001 & 1.55e-014 & 4.9933773e-001 & 1.62e-14 \\ \hline
 1.9939880e-001 & 5.90e-014 & 5.0082218e-001 & 1.04e-12 \\ \hline
 1.9813769e-001 & 1.08e-009 & 4.9764898e-001 & 8.76e-11 \\ \hline
 \multicolumn{2} {|c|} {$\tau=0.6$} & \multicolumn{2} {|c|} {$\tau=0.8$}
 \\ \hline
 $\rho_j$ & $|\rho_j-\sigma_j|/\sigma_j$ & $\rho_j$ & $|\rho_j-\sigma_j|/\sigma_j$ \\ \hline
   6.0106012e-001 & 4.29e-12 & 8.0014466e-001 & 2.41e-12 \\ \hline
   6.0193472e-001 & 4.22e-11 & 7.9954438e-001 & 5.46e-12 \\ \hline
   5.9689466e-001 & 2.40e-13 & 7.9932106e-001 & 1.08e-10 \\ \hline
\end{tabular}
\end{table}

\begin{table}[htp]
\caption{DW2048 for $k=3$, $m=30, 40, 50$, $\tau=0.2, 0.5, 0.6,
0.8$} \label{tab3}
\begin{tabular}{|c|c|c|c|c|c|c|c|c|} \hline
&\multicolumn{4} {|c|} {$\tau=0.2$}&\multicolumn{4} {|c|}
{$\tau=0.5$} \\ \hline
$m$ & $iter$ & $time$ & $mv$ & $stopcrit$ &
$iter$ & $time$ & $mv$ & $stopcrit$ \\ \hline
 30 & 501 & 109 & 11655 & 9.97e-007 & 255 & 51.2 & 6123 & 9.93e-007 \\ \hline
 40 & 298 & 113 & 9993 & 9.81e-007 & 97 & 38.2 & 3271 & 9.90e-007 \\ \hline
 50 & 221 & 136 & 9652 & 9.85e-007 & 83 & 50.9 & 3656 & 9.10e-007 \\ \hline
&\multicolumn{4} {|c|} {$\tau=0.6$}&\multicolumn{4} {|c|}
{$\tau=0.8$}\\ \hline $m$ & $iter$ & $time$ & $mv$ & $stopcrit$ &
$iter$ & $time$ & $mv$ & $stopcrit$ \\ \hline
 30 & 125 & 24.6 & 3006 & 9.11e-007 & 405 & 78.2 & 9525 & 9.87e-007 \\ \hline
 40 & 69 & 27.0 & 2343 & 9.47e-007 & 180 & 70.4 & 6079 & 9.64e-007 \\ \hline
 50 & 46 & 28.1 & 2012 & 8.61e-007 & 136 & 81.8 & 5989 & 9.50e-007 \\ \hline
\end{tabular}
\end{table}

\subsection{Computation of interior singular triplets for different $k$}

We compute $k=1,3,5,10$ smallest singular triplets nearest $\tau=4.5$
of LSHP2233, a $2233 \times 2233$ matrix. Table \ref{tab5} reports
the results. We see that our algorithm can compute the desired
singular triplets with high precision.

\begin{table}[htp]
\caption{Ten computed singular values of LSHP2233 nearest $\tau=4.5$
for $m=50$} \label{tab4}
\begin{tabular}{|c|c|c|c|} \hline
 $\rho_1$  & $|\rho_1-\sigma_1|/\sigma_1$ & $\rho_2$ &
 $|\rho_2-\sigma_2|/\sigma_2$ \\ \hline
 4.4988631 & 1.58e-15 & 4.5091282 & 1.36e-14 \\ \hline
 $\rho_3$  & $|\rho_3-\sigma_3|/\sigma_3$ & $\rho_4$ &
 $|\rho_4-\sigma_4|/\sigma_4$ \\ \hline
 4.5113859  & 1.22e-14 & 4.4815289 & 6.54e-15 \\ \hline
 $\rho_5$  & $|\rho_5-\sigma_5|/\sigma_5$ & $\rho_6$ &
 $|\rho_6-\sigma_6|/\sigma_6$ \\ \hline
 4.5188882  & 5.11e-15 & 4.5210494 & 1.18e-14 \\ \hline
 $\rho_7$  & $|\rho_7-\sigma_7|/\sigma_7$ & $\rho_8$ &
 $|\rho_8-\sigma_8|/\sigma_8$ \\ \hline
 4.4783693  & 1.07e-14 & 4.4716358 & 8.74e-15 \\ \hline
 $\rho_9$  & $|\rho_9-\sigma_9|/\sigma_9$ & $\rho_{10}$ &
 $|\rho_{10}-\sigma_{10}|/\sigma_{10}$ \\ \hline
 4.5331457  & 5.68e-15 & 4.4638926 & 2.03e-10 \\ \hline
\end{tabular}
\end{table}

\begin{table}[htp]
\caption{LSHP2233 for $k=1, 3, 5, 10$, $m=30, 40 , 50$, $\tau=4.5$}
\label{tab5}
\begin{tabular}{|c|c|c|c|c|c|c|c|c|} \hline
&\multicolumn{4} {|c|} {$k=1$}&\multicolumn{4} {|c|} {$k=3$} \\ \hline
$m$ & $iter$ & $time$ & $mv$ & $stopcrit$ & $iter$ & $time$ & $mv$ & $stopcrit$ \\ \hline
 30 & 467 & 108 & 11920 & 6.88e-006 & 560 & 127 & 13171 & 6.95e-006 \\ \hline
 40 & 190 & 83.6 & 6844 & 6.86e-006 & 230 & 98.1 & 7826 & 6.79e-006 \\ \hline
 50 & 159 & 114 & 7188 & 6.63e-006 & 216 & 140 & 9404 & 6.66e-006 \\ \hline
&\multicolumn{4} {|c|} {$k=5$}&\multicolumn{4} {|c|} {$k=10$} \\ \hline
$m$ & $iter$ & $time$ & $mv$ & $stopcrit$ & $iter$ & $time$ & $mv$ & $stopcrit$ \\ \hline
 30 & 322 & 64.8 & 6972 & 6.99e-006 & 651 & 103 & 10761 & 6.91e-006 \\ \hline
 40 & 207 & 78.5 & 6632 & 6.78e-006 & 168 & 61.0 & 4548 & 4.70e-006 \\ \hline
 50 & 132 & 83.8 & 5487 & 6.40e-006 & 165 & 92.4 & 6003 & 6.99e-006 \\ \hline
\end{tabular}
\end{table}

\section{Conclusion}

In this paper, combining the harmonic projection principle with the
implicit restarting technique, we propose an implicitly restarted
harmonic Lanczos bidiagonalization algorithm for computing some
interior singular triplets. Based on Morgan's harmonic shift strategy for computing interior eigenpairs, 
we give a selection of the shifts within our algorithm.
Numerical experiments show that our algorithm is suitable for interior SVD problems. 
The interior singular values can be computed with higher relative precision. 

The Matlab code can be obtained from the authors upon request.

\begin{acknowledgements}

We thank Baglama and Reichel very much for generously providing
their Matlab code of Lanczos Bidiagonalization process on their
homepage, which reduces our programming work greatly.

\end{acknowledgements}




\end{document}